\documentclass[12pt]{amsart}
\usepackage{amssymb}
\usepackage{amsmath}
\usepackage{amsfonts}

\newtheorem{theorem}{Theorem}

\newtheorem{lemma}[theorem]{Lemma}

\newcommand\C{\mathbb{C}}

\newcommand\R{\mathbb{R}}
\newcommand\Z{\mathbb{Z}}
\newcommand\wrt{{\ d}}

\renewcommand\Re{\operatorname{Re}}

\begin{document}

\title[Bochner-Riesz summability for ...]{Bochner-Riesz summability for analytic functions on the 
$m$-complex unit sphere and for cylindrically symmetric functions 
on $\R^{n-1} \times \R$}

\author{Adam Sikora}
\address{Adam  Sikora, Department of Mathematical Sciences, 
New Mexico State University, Las Cruces, NM 88003-8001, USA}
\email{asikora@nmsu.edu}

\author{Terence Tao}
\address{Terence Tao, Department of Mathematics, UCLA, LA CA 90095 USA}
\email{tao@math.ucla.edu}


\begin{abstract}
We prove that spectral projections of Laplace-Beltrami operator on the 
$m$-complex unit sphere $E_{\Delta_{S^{2m-1}}}([0,R))$ are 
uniformly bounded as an
operator from $H^p(S^{2m-1})$ to $L^p(S^{2m-1})$ for all $p\in
(1,\infty)$.
We also show that the Bochner-Riesz conjecture is true when restricted
to cylindrically symmetric functions on $\R^{n-1} \times \R$.
\end{abstract}

\thanks{This work was conducted at ANU. \hfill{} \linebreak 
AS is partially  supported by a Minigrant and 
Summer Research Award  from New Mexico State
University.
\hfill{} \linebreak   
TT is a Clay Prize Fellow and is supported by a grant from the Packard Foundation.  
} 

\maketitle

Suppose that $L$ is a positive definite, self-adjoint operator acting
on $L^2(X,\mu)$, where $X$ is a measurable space with a measure
$\mu$. Such operator admits a spectral resolution 
$$
L= \int_0^\infty \lambda \wrt E_L(\lambda).
$$
By the spectral theorem, if $F$ is a Borel bounded function on $[0,\infty)$,
then the operator $F(L)$ given by 
\begin{equation}\label{equw}
F(L) := \int_0^\infty F(\lambda) \wrt E_L(\lambda)
\end{equation}
is well-defined and bounded on $L^2(X,\mu)$.
One of the fundamental problems in the theory of spectral multipliers is to 
determine when $F(L)$ is bounded on $L^p$ for some $p\neq 2$.  
An interesting example is the following family of functions
\begin{equation}\label{vab1}
   S^{\delta}_R(\lambda) :=
      \left\{
       \begin{array}{cl}
       (1-\lambda/R)^{\delta}  &\mbox{for}\;\; \lambda \le R \\
       0  &\mbox{for}\;\; \lambda > R. \\
       \end{array}
      \right.
   \end{equation}
We  define the operator $S^{\delta}_R(L)$ using (\ref{equw}).
$S^{\delta}_R(L)$ is called the Riesz mean or the Bochner-Riesz mean of
order $\delta$. The basic question in the theory of the Riesz means is  to
establish  the critical exponent for the continuity and convergence of
the Riesz means. More precisely we want to study the optimal range of
$\delta$ for which the Riesz means $S^{\delta}_R(L)$ are uniformly
bounded on $L^p(X)$, or in other words that
\begin{equation}\label{unn}
\sup_{R>0}\|S^\delta_R(L)\|_{L^p(X,\mu)\to L^p(X,\mu)} <\infty.
\end{equation}
Since the publication
of Riesz's paper \cite{Ri} the summability  of  the  Riesz means has been
one of the most fundamental problems in Harmonic Analysis, see e.g.
\cite[IX.2 and \S IX.6B]{St2}.
The best understood case is $X=\R^n$ with the Lebesgue measure
and $L=\Delta_{\R^n}$, where $\Delta_{\R^n}$ is the standard Laplace operator.
It is known that $S^0_R(\Delta_{\R^1})$ is uniformly bounded on $L^p$
for all $1<p<\infty$ but if $n \neq 1$, then 
$S^0_R(\Delta_{\R^n})$ is bounded on $L^p(\R^n)$ only when $p=2$, see
\cite{Fe}.

Next denote by $L^p_{rad}(\R^n)$ the subspace of elements of
$L^p(\R^n)$ which are invariant under the action on $\mbox{\rm SO}(n)$.
In \cite{He} Herz proved 
\begin{theorem}\label{herz}
Suppose that $(n-1)/2 \ge \delta \ge 0$. Then
$$
\sup_{R>0}\|S^\delta_R(\Delta_{\R^N})\|_{L^p_{rad}(X,\mu)\to L^p(X,\mu)} <\infty
$$
if and only if 
\begin{equation}\label{p}
\frac{2n}{n+1+2\delta} < p < \frac{2n}{n-1-2\delta}.
\end{equation}
\end{theorem}
Thus if $S^\delta_R(\Delta_{\R^n})$ 
is uniformly continuous on $L^p(\R^n)$, then $p$ satisfies (\ref{p}).
It is conjectured that for $\delta >0$ the converse is true, 
but this is only known
for $n=2$.
\begin{theorem}\label{rm}
Suppose that $L =\Delta_{\R^n}$ or $L=\Delta_{S^n}$ and $X=\R^n$ or
$X=S^n$ respectively. Next 
suppose that $(n-1)/2 \ge \delta > 0$, $p$ satisfies (\ref{p})
and that 
\begin{equation}\label{p1}
1 \le  p \le  \frac{2(n+1)}{n+3}\quad \mbox{or} \quad
\frac{2(n+1)}{n-1}\le p \le \infty.
\end{equation}
Then
\begin{equation}\label{unr}
\sup_{R>0}\|S^\delta_R(L)\|_{L^p(X,\mu)\to L^p(X,\mu)} <\infty .
\end{equation}
If $n=2$ then (\ref{unr}) holds if and only if $p$ satisfies (\ref{p}).
\end{theorem}
For a proof see \cite[Chapter IX]{St2} or  \cite[Chapter 5.2]{So}.
It is an open problem if condition (\ref{p1}) can be removed from 
the assumptions of the Theorem~\ref{rm} also for $n>2$. This problem
has been studied
intensively  and will probably
remain open in the near future. Therefore it seems to be natural to
investigate simpler versions of problem, where we consider a closed subspace 
$Q \subset L^p(X,\mu)$ and ask if the operators $S^\delta_R(L)$ are
uniformly continuous from $Q$ to $L^p(x,\mu)$. If we consider 
the subspace $Q=L^p_{rad}(\R^n)$, then a complete answer is provided by 
Theorem~\ref{herz}. Using standard techniques it is easy to prove
the same result for  $Q=L^p_{rad}(S^n)$.

Here we propose to study the following examples:

\begin{itemize}

\item The subspace $H^p(S^{2m-1})\subset L^p(S^{2m-1})$ of analytic functions on the $m$-complex sphere $S^{2m-1}$.

\item Cylindrically symmetric functions on $\R^{n-1} \times \R$.

\end{itemize}

In both cases we obtain sharp results.  In the first case this will be because the Laplacian on the sphere can be rewritten in terms of a single vector field when restricted to analytic functions; in the second case we shall use a Radon transform to reduce the problem to a two-dimensional weighted Bochner-Riesz estimate, together with some one-dimensional weighted estimates which follow from Hardy's inequality.

\subsection*{Analytic functions on $m$-complex sphere} 

To be more precise we have to introduce more notation.
We use $z$ as a general element of $\C^m$, i.e. $z=(z_1,\ldots,z_m)$,
$z_i\in \C$, $i=1,2,\ldots,m$, $m \ge 2$. We let $\bar{z}:=(\bar{z}_1,
\ldots,\bar{z}_m)$ denote the conjugate of $z$, and
and  $B_m :=\{z\in \C^m \colon |z| \le 1\}$ denote the unit ball, where $|z|:=\big(\sum_{i=1}^{m}|z_i|^2\bar)^{1/2}$.
We use  $S^{2m-1} := \partial B_m := \{z\in \C^m \colon |z|= 1\}$ to 
denote the  boundary of $B_m$,

Let $\wrt \sigma$ be the Lebesgue area element of the unit sphere $S^{2m-1}$.
By $A(\partial B_m)$ we denote the class of all
$f\in C(\partial B_m)$ that are restrictions to $\partial B_m$ of an analytic
function on $B$.
For $ 0<p< \infty$ we let $H^p(S^{2m-1})$ be the $L^p$-closure of $A(\partial
B_m)$, see \cite[Definition 5.6.7]{Ru}.
Our first result is somehow surprising given Theorems~\ref{herz}~and~\ref{rm}.
\begin{theorem}\label{unt}
For any natural number $m$ and for any  $1< p < \infty$ 
\begin{equation}\label{un}
\sup_{R>0}\|S^0_R(\Delta_{S^{2m-1}})\|_{H^p(S^{2m-1})\to L^p(S^{2m-1})}
    < \infty.
\end{equation}
Moreover the operators $S^0_R(\Delta_{S^{2m-1}})$ are uniformly
continuous from $H^1(S^{2m-1})$ to $L^{1,\infty}$ (weak-$L^1$). This means that\footnote{As
usual we will use $C$ to denote various constants which depend on parameters
such as the ambient dimension $m$ and exponents such as $p$.}
$$
\sigma(\{ x \in S^{2m-1} : |S^0_R(\Delta_{S^{2m-1}}) f(x)|)  > \lambda \}
\leq C \, \frac{\|f \|_{H^1}}{\lambda}
\qquad\forall \lambda, R \in \R^+ \quad\forall f \in H^1(S^{2m-1}).
$$
\end{theorem}

Before we prove Theorem~\ref{unt} we want to discuss some results
described in \cite{CQ}. As we will see these results are closely related
to Theorem~\ref{unt}.

Let $k$ be a nonnegative integer and let $\nu=(k_1,\ldots,k_m)$,
where $k_1+\ldots +k_m=k$ and $k_j \ge 0$. We define polynomials 
$p_\nu^k$ by the formula
$$
p_\nu^k := \Big(\frac{k!}{k_1!\ldots k_m!}\Big)^{1/2}z_1^{k_1}\ldots z_m^{k_m}.
$$
We now recall some standard facts from \cite{CQ,Ru}.  First, we have the following
orthonormal basis property.
The following theorem is well known:
\begin{theorem}
The system of functions $\{p_{\nu}^k\}$ is an orthonormal (but not complete)
system in $L^2(\partial B_m)$. The system $\{p_{\nu}^k\}$ is orthonormal
and complete in $H^2(\partial B_m)$.
\end{theorem}
Also, we have the reproducing formula
$$
\sum_{|\nu|=k}p_{\nu}^k(z) \overline{p_{\nu}^k}(\xi)=\frac{1}{\omega_{2m-1}}
\binom{m+k-1}{k}(z\bar{\xi})^k,
$$
where $z\bar{\xi}=\sum_{k=1}^{m}z_k\bar{\xi}_k$ and $\omega_{n-1}=
2\pi^{n}/\Gamma(n/2)=\sigma(S^{n-1})$.
For given function $b\colon \R_+ \to \C$ we define the multiplier operator 
$M_b$ by the formula
$$
M_bf(z) :=\int_{\partial B_m} \sum_{k=0}^{\infty}b(k)
\sum_{|\nu|=k}p_{\nu}^k(z) \overline{p_{\nu}^k}(\xi) f(\xi)=
\int_{\partial B_m} 
\sum_{k=0}^{\infty}b(k)f(\xi)\binom{m+k-1}{k}(z\bar{\xi})^k.
$$ 
Let $f\in L^p(\partial B_m)$ , $1\le p \le \infty$. The Cauchy
integral $C(f)$ of $f$
$$
C(f)(z) :=\int_{\partial B_m}\frac{f(\xi)}{(1-z\bar{\xi})^m} \wrt \sigma(\xi)
$$
is then well defined and holomorphic in $B_m$. The operator 
$$
P(f)(\xi) := \lim_{r \to 1_-}C(f)(r\xi)
$$
is the projection of $L^p(\partial B_m)$ onto the Hardy space 
$H^p(\partial B_m)$ and is bounded from $L^p(\partial B_m)$ to 
$H^p(\partial B_m)$, $1 <p <\infty$. 
Note that if we put $b=1$, then $P=M_b$. 
Now if we define the radial
Dirac operator $D$ acting on $H^2(\partial B_m)$ by the formula
$$D := \sum_{k=1}^{m} z_k\frac{\partial}{\partial z_k},$$
then the operator $DP$ is a self-adjoint positive defined operator 
acting on $L^2(\partial B_m)$ and $M_b= b(DP)$, where $b(DP)$ is
defined by (\ref{equw}).

We define the Banach space $H^\infty(\Sigma_\omega)$ by
$$
H^\infty(\Sigma_\omega) := \{ b \colon \Sigma_\omega \to \C \quad |\quad  b \,\,
\mbox{is holomorphic in} \,
\Sigma_\omega \, \mbox{and} \, \|f\|_{L^\infty(\Sigma_\omega)} < \infty \},
$$
where $\Sigma_\omega$ is the sector
$$\Sigma_\omega :=
\{z \in \C \quad |\quad  |\mbox{\rm arg} \, z|
\le \omega \}.$$
The following theorem is the main result obtained in \cite{CQ}:
\begin{theorem}\cite{CQ}
Suppose that $b \in H^\infty(\Sigma_\omega)$. Then 
the operator $M_b$ can be extended to a bounded operator from
$L^p(\partial B_m)$ to $L^p(\partial B_m)$, $1<p <\infty,$  and from
$L^1(\partial B_m)$ to $L^{1,\infty}(\partial B_m)$.
\end{theorem}
Finally we  state a H\"ormander type multiplier
theorem for the operator $DP$, see \cite[Theorem~7.9.5, pp. 243]{Her4}.
\begin{theorem}\label{her}
Suppose that for some exponent $s>1/2$
function $b$ satisfies the H\"ormander type
condition
\begin{equation}\label{herr}
\sup_{ t  > 0}  \Vert  \eta \, \delta_t b \Vert_{H_{s}} < \infty,
\end{equation}
where $ \delta_t b(\lambda)=b(t\lambda)$,
$H_s$ is the Sobolev space of order $s$ and $\eta \in C_c^{\infty}(\R_+)$ is a fixed  function,  not  identically
   zero. Then 
the operator $M_b$ can be extended to a bounded operator from
$L^p(\partial B_m)$ to $L^p(\partial B_m)$, $1<p <\infty$.
\end{theorem}
By the  Cauchy formula 
\begin{equation}\label{3.1e}
\sup_{ t  > 0} \|\eta \,  \delta_t b \|_{H_{k}}\le C\sup_{\lambda > 0}  \vert  \lambda^k b^{(k)}(\lambda) \vert  \le
\frac{C_k}{\omega^k}
\Vert b \Vert_{L^\infty(\Sigma_\omega)}, \ \forall k \in \Z_+.
\end{equation}
Hence by (\ref{3.1e}) and interpolation for all $s>0$
\begin{equation}\label{3.6e}
\sup_{ t  > 0}  \Vert  \eta \, \delta_t b \Vert_{H_{s}}  \le
\frac{C_{s}}{\omega^{s}}\Vert b \Vert_{L^\infty(\Sigma_\omega)}.
\end{equation}
>From (\ref{3.6e}) we see that  
Theorem~\ref{her} strengthens the result of Cowling and Qian in a significant
way. Note that the critical exponent $1/2$ in Theorem~\ref{her} is optimal.
Theorem~\ref{unt} and Theorem~\ref{her}  have a very similar proof so
we prove them simultaneously\footnote{Essentially Theorem~\ref{unt} follows from Theorem~\ref{her}}.

\begin{proof}[Proof of Theorem~\ref{unt} and Theorem~\ref{her}]
Let us recall the well known Bochner formula
\begin{equation}\label{ww}
\int_{\partial B_m}f \wrt\sigma=\int_{\partial B_m}\wrt\sigma(\xi)
\frac{1}{2\pi}\int_\pi^\pi f(e^{it}\xi)\wrt t,
\end{equation}
see \cite[Proposition 1.4.7]{Ru}. If $f \colon \partial B_m \to \C$,
then for all $\xi \in \partial B_m$ 
we  define $f_\xi \colon \partial B_1 \to \C$ by the formula
$$
f_\xi(t):=f(e^{it}\xi).
$$
Now we can rewrite (\ref{ww}) in the following way
\begin{equation}\label{www}
\|f\|_{L^p(S^{2m-1})}^p=\int_{\partial B_m} \wrt \sigma(\xi)
\|f_\xi\|_{L^p(S^{1})}^p
\end{equation}
Next we define the self-adjoint operator $iZ_m$ by the formula
$$
iZ_mf(z):=i\frac{\wrt}{ \wrt t}f(ze^{it})|_{t=0}.
$$ 
Note that for any $f_\xi \colon \partial B_1 \to \C$
$$
\exp(-itZ_m)f(z)=f(e^{it}z).
$$
Also, for any polynomial $p^k_\nu$ one can verify the identity
$$ \Delta_{S^{2m-1}} p^k_\nu = k (k+m) p^k_\nu = iZ_m (iZ_m + m) p^k_\nu$$
and hence that 
$$
S^0_R(\Delta_{S^{2m-1}})P=S^0_{R(R+m)}(iZ_m)P=\chi_{[0,R(R+m))}(DP)
$$
and 
$$
b(DP)=M_b=b(iZ_m)P.
$$
Moreover if $g=b(iZ_m)f$ then $g_\xi=b(iZ_1)f_\xi$ and
\begin{eqnarray*}
\|b(DP)f\|_{L^p(S^{2m-1})}^p=\int_{\partial B_m} \wrt \sigma(\xi)
\|b(iZ_1)(Pf)_\xi\|_{L^p(S^{1})}^p \\ \le 
\|b(iZ_1)\|_{L^p(S^{1}) \to L^p(S^{1}) }^p \int_{\partial B_m} \wrt \sigma(\xi)
\|(Pf)_\xi\|_{L^p(S^{1})}^p\\=\|b(iZ_1)\|_{L^p(S^{1}) \to L^p(S^{1}) }^p
\|Pf\|_{L^p(S^{2m-1})}^p
\end{eqnarray*}
To end the proof of Theorem~\ref{her} we note that by the H\"ormdander 
multiplier theorem
$\|b(iZ_1)\|_{L^p(S^{1})\to L^p(S^{1}) }   \le C
\sup_{ t  > 0}  \Vert  \eta \, \delta_t b \Vert_{H_{s}}$,
 see \cite[Theorem~7.9.5, pp. 243]{Her4}. To prove (\ref{un}) we note
that $\sup_R\|S^0_R(iZ_1)\|_{L^p(S^{1})\to L^p(S^{1}) }<\infty$ by 
Theorem~\ref{herz}.

Finally note that 
$$
\sigma(\{ z\in \partial B_m \colon f(z)| > \lambda \})=
\int_{\partial B_m} \wrt \sigma(\xi)
 \sigma_1(\{ e^{it} \in S^1 \colon   |f_\xi(t)|> \lambda \}). 
$$
Hence if $f \in H^1(S^{2m-1})$, then
$$
\sigma(\{ z\in \partial B_m \colon |b(DP) f(z)| > \lambda \})
 \le \frac{1}{\lambda} \|b(iZ_1)\|_{L^1(S^{1}) \to L^{1,\infty}(S^{1})}.
$$  
\end{proof}
{\em Remark.} We proved that the norm of operator $b(DP)$ from 
$H^1(S^{2m-1})$ to $L^{1,\infty}(S^{2m-1})$ 
is bounded by the weak $(1,1)$ norm of
the operator $b(iZ_1)$. One can ask if the operator $b(DP)$
is bounded from $L^1(S^{2m-1})$ to $L^{1,\infty}(S^{2m-1})$.
The result of Cowling and Qian says that this is the case 
for  $b \in H^\infty(\Sigma_\omega)$. However, this is not usually
the case for the functions $b$ which are only assumed to satisfy (\ref{herr}).

\subsection*{Cylindrically symmetric functions on $\R^{n-1} \times \R$}

Let $n \geq 3$.  We now work in the space $\R^{n-1} \times \R := \{ (x,t): x \in \R^{n-1}, t \in \R \}$, and consider the rotation group $SO(n-1)$ acting on the first factor $\R^{n-1}$.  
Let $L^p_{cyl}(\R^{n-1} \times \R)$ denote those
functions $f$ in $L^p$ which are invariant with respect to the $SO(n-1)$
action (and are therefore cylindrically symmetric).  Observe that the Laplacian
$\Delta_{\R^{n-1} \times \R}$, and hence all spectral multipliers based on
this Laplacian, commute with $SO(n-1)$ and hence preserve the space of
cylindrically symmetric functions.

We now show that the Bochner-Riesz conjecture is true when restricted
to cylindrically symmetric functions on $\R^{n-1} \times \R$.  In other
words we show that
$$ 
\sup_R \| S^\delta_R( \Delta_{\R^{n-1} \times \R} ) \|_{L^p_{cyl}(\R^{n-1}%
  \times \R) \to L^p_{cyl}(\R^{n-1} \times \R)} < \infty 
$$
whenever (4) holds.

\begin{proof}
By scaling we can take $R=1$.  The claim is known to be true when
$p=1$, $p=2$ or $p=\infty$; by duality and interpolation it then
suffices to show that
\begin{equation}\label{dd0}
\| S^\delta_1( \Delta_{\R^{n-1} \times \R} ) \|_{L^{2n/(n+1)}_{cyl}(\R^3 \times
  \R) \to L^{2n/(n+1)}_{cyl}(\R^3 \times \R)} < \infty
\end{equation}
for any $\delta > 0$.

Fix $\delta$; our constants $C$ may depend on $\delta$.  By duality it suffices to show that
\begin{equation}\label{dd}
|\langle S^\delta_1 ( \Delta_{\R^{n-1} \times \R} ) f, g \rangle|
\lesssim \| f \|_{L^{2n/(n+1)}_{cyl}(\R^{n-1} \times \R)} \| g \|_{L^{2n/(n-1)}_{cyl}(\R^{n-1} \times \R)}
\end{equation}
for all cylindrically symmetric test functions $f, g$.  Without loss of generality we may assume that $f$, $g$ are real.

We now use the method of descent to reduce this problem from $n$ dimensions to $2$, exploiting the cylindrical symmetry.  The method here (inspired by work of Rubio de Francia) is not specific to the Bochner-Riesz multipliers and could be extended to other cylindrically symmetric multipliers (and to other types of cylindrical symmetry).

We introduce the $n$-dimensional Fourier transform
$$ {\mathcal F}_n f(\xi, \tau) := \int_\R \int_{\R^{n-1}} e^{-2\pi i (x \cdot \xi + t \tau)} f(x, t)\wrt x \wrt t$$
for $\xi \in \R^{n-1}$, $\tau \in \R$, and observe that ${\mathcal F}_n$ is cylindrically symmetric in the $\xi$ variable if $f$ is cylindrically symmetric in the $x$ variable.  By Plancherel, we can write the left-hand side of \eqref{dd} as
$$ |\int_\R \int_{\R^{n-1}} (1-4\pi^2|(\xi,\tau)|)_+^\delta {\mathcal F}_n f(\xi,\tau) \overline{{\mathcal F}_n g(\xi,\tau)}\wrt \xi \wrt \tau|.$$
We parameterize $\R^{n-1}$ as $(x_1, \underline x)$, with $x_1 \in \R$, $\underline x \in \R^{n-2}$, and write $e_1 := (1,0)$.  By a change to polar co-ordinates and exploiting the cylindrical symmetry of the integrand, we may write the previous expression as
\begin{equation}\label{pr}
\omega_{n-2} |\int_\R \int_{\R} (1-|(\xi_1 e_1,\tau)|)_+^\delta {\mathcal F}_n f(\xi_1 e_1,\tau) \overline{{\mathcal F}_n g(\xi_1 e_1,\tau)}\ |\xi_1|^{n-2} d\xi_1 d\tau|.
\end{equation}
Introduce the projection operator $P$ defined by
$$ Pf(x_1, t) := \int_{\R^{n-2}} f((x_1, \underline x), t)\wrt \underline{x}$$
and the two-dimensional Fourier transform
$$ {\mathcal F}_2 f(\xi_1, \tau) := \int_\R \int_{\R} e^{-2\pi i (x_1 \xi_1 + t \tau)} f(x_1, t)\wrt x_1 \wrt t,$$
where we have parameterized $\R \times \R$ by $\{ (x_1,t): x_1 \in \R, t \in \R \}$.  Observe the identity
$$ {\mathcal F}_n f(\xi_1 e_1, \tau) = {\mathcal F}_2 Pf(\xi_1, \tau).$$
Thus we may write (\ref{pr}) as
$$ \omega_{n-2} |\int_\R \int_{\R} (1-4\pi^2|(\xi_1 e_1,\tau)|)_+^\delta {\mathcal F}_2 Pf(\xi_1 e_1,\tau) \overline{{\mathcal F}_2 Pg(\xi_1,\tau)}\ |\xi_1|^{n-2} d\xi_1 d\tau|.$$
Distributing the ``derivatives'' $|\xi_1|^{n-2}$ equally among $Pf$ and $Pg$, we can rewrite the previous as
$$ \omega_{n-2} |\int_\R \int_{\R} (1-4\pi^2|(\xi_1 e_1,\tau)|)_+^\delta {\mathcal F}_2 |\partial_{x_1}|^{(n-2)/2} Pf(\xi_1 e_1,\tau) \overline{{\mathcal F}_2 |\partial_{x_1}|^{(n-2)/2} Pg(\xi_1,\tau)}\wrt \xi_1 \wrt \tau|,$$
which after undoing the Plancherel becomes
$$ \omega_{n-2} | \langle S^\delta_1 ( \Delta_{\R \times \R} ) |\partial_{x_1}|^{(n-2)/2} Pf, |\partial_{x_1}|^{(n-2)/2} Pg \rangle|.$$

We claim the three weighted estimates\footnote{The exponents here may seem somewhat arbitrary, 
but they are forced on us by scaling considerations.  Note that the number of derivatives $(n-2)/2$ 
in \eqref{l1}, \eqref{l2} matches the amount of smoothing available for the hyperplane Radon transform 
in $\R^{n-1}$, see e.g. \cite{St2}; this Radon transform is essentially equivalent to $P$ for 
cylindrically symmetric functions.  This argument was inspired by a simpler version which was available 
in the four-dimensional case, which we describe in the next section.}
\begin{align}
\| |x_1|^{\frac{n-2}{2n}} |\partial_{x_1}|^{(n-2)/2} Pf \|_{L^{2n/(n+1)}(\R \times \R)}
&\leq C \| f \|_{L^{2n/(n+1)}_{cyl}(\R^{n-1} \times \R)} \label{l1} \\
\| |x_1|^{-\frac{n-2}{2n}} |\partial_{x_1}|^{(n-2)/2} Pg \|_{L^{2n/(n-1)}(\R \times \R)}
&\leq C \| f \|_{L^{2n/(n-1)}_{cyl}(\R^{n-1} \times \R)} \label{l2} \\
\| |x_1|^{\frac{n-2}{2n}} S^\delta_1 ( \Delta_{\R \times \R} ) |x_1|^{-\frac{n-2}{2n} } \|_{L^{2n/(n+1)}(\R \times \R) \to L^{2n/(n+1)}(\R \times \R)} &\leq C;\label{l3}
\end{align}
the claim \eqref{dd} then clearly follows from composing the above estimates and applying H\"older's inequality.

It remains to prove \eqref{l1}, \eqref{l2}, \eqref{l3}.  We begin with \eqref{l3}.  We shall in fact prove the more general

\begin{lemma}\label{delta}  For any $\delta > 0$ and $4/3 \leq p < 2$, the operator
$|x_1|^{\alpha} S^\delta_1( \Delta_{\R \times \R} ) |x_1|^{-\alpha}$ is bounded in $L^p( \R \times \R )$ whenever
$$ |\alpha| \leq \frac{3}{2} - \frac{2}{p}.$$
\end{lemma}

Clearly \eqref{l3} follows from this lemma by setting $p := 2n/(n+1)$ and
$\alpha := (n-2)/2n = 3/2 - 2/p$.

\begin{proof}
Since the Bochner-Riesz conjecture is known in two dimensions, 
see \cite[Chapter IX]{St2} or \cite[Chapter 5.2]{So}, 
the claim is true when $p=4/3$, $\alpha = 0$, and $\delta > 0$ is arbitrary.

Next, we assert that the claim is true when 
$1 \leq p \le 2$, $ -1/p <-\alpha  <1- 1/p$, and $\delta > 100$.
  To see this, observe that the convolution kernel $K^\delta_1$ of 
$S^\delta_1( \Delta_{\R \times \R})$ decays rapidly, say 
$$
|K^\delta_1(x)| \le C( (1 + |x|)^{-20} )\le C(1+|x_1|)^{-10}(1+|x_2|)^{-10}
.$$  
Hence it suffices to show that the operator 
$P_\alpha \colon \; f(x) \mapsto |x|^\alpha \int_{\R}C(1+|x-y|)^{-10} |y|^{-\alpha} f(y)dy$ 
is bounded on $L^p(\R)$ for all  $ -1/p <-\alpha  <1- 1/p$. 
But this follows from \cite[Corollary p. 205 and \S 6.4 p. 218]{St2}.
Alternatively to show that $P_\alpha $ 
is bounded on $L^p(\R)$  for all $ -1/p < -\alpha  <1- 1/p$ we note that it
is true for $p=1$ and $p=\infty$ because
$$
\sup_x |x|^\alpha \int_{\R}C(1+|x-y|)^{-10} |y|^{-\alpha} \le C <\infty
$$
for $0<\alpha <1$. Then we get the boundedness of $P_\alpha$
  for all $ -1/p <-\alpha  <1- 1/p$
by interpolation.

>From the previous observations and complex interpolation\footnote{More precisely, by interpolating the $p=2$ estimates with the $p=4/3$ ones we obtain the theorem assuming that $|\alpha|$ is strictly less than $3/2-2/p$.  To obtain the
endpoint case (which is what we need for \eqref{l3}) for we exploit the fact that $\delta > 0$ and interpolate with the
$\delta > 100$ estimates; note that the range $-1/p < -\alpha < 1-1/p$
is strictly larger than $|\alpha| \leq 3/2 - 2/p$.}
 we see that it suffices to prove the claim when $p=2$, $\delta = 0$, and $|\alpha| < 1/2$.  In other words, we need to show
$$ \| |x_1|^\alpha S^0_1( \Delta_{\R \times \R} ) |x_1|^{-\alpha} f \|_{L^2(\R \times \R)}
\lesssim \| f \|_{L^2(\R \times \R)}.$$
We write $f(x_1,t)$ using a Fourier transform in time as
$$ f(x_1,t) = \int g_\tau(x_1) e^{i t \tau}\wrt \tau$$
for some function $g_\tau$.  Observe the identity
$$ S^0_1( \Delta_{\R \times \R} ) |x_1|^{-\alpha} f(x_1,t) = \int_{|\tau| \leq 1}
   S^0_{1-\tau^2}(\Delta_\R) |x_1|^{-\alpha} g_\tau(x_1) e^{i t \tau}\wrt \tau$$
(just take the Fourier transforms of both sides).  From a Plancherel in
the $t$ variable we thus reduce to showing the one-dimensional weighted estimate
$$ 
\| |x_1|^\alpha S^0_{1-\tau^2}(\Delta_\R) |x_1|^{-\alpha} g_\tau \|_{L^2(\R)}
\leq C \| g_\tau \|_{L^2(\R)}
$$
uniformly in $\tau$.  But the operator $S^0_{1-\tau^2}(\Delta_\R)$ is
just a linear combination of modulated Hilbert transforms, and
$|x_1|^{2\alpha}$ is an $A_2$ weight for $-1/2 < \alpha < 1/2$, so the claim follows, see \cite[Corollary p. 205 and \S 6.4 p. 218]{St2}.
\end{proof}

It remains to prove \eqref{l1}, \eqref{l2}.  We remark that in the special case $n=4$ these identities are easy (especially if we replace $|\partial_{x_1}|$ by $\partial_{x_1}$; we will return to this point in the next section).  For general dimension $n$, we begin by observing from use of polar co-ordinates in $\R^{n-2}$ that
$$ Pf(x_1, t) = \omega_{n-3} \int_0^\infty f(\sqrt{x_1^2 + r^2}, t) r^{n-3}\wrt r$$
for cylindrically symmetric $f$, where we have abused notation and written $f(|x|, t)$ for $f(x, t)$.  Making the change of variables $y := \sqrt{x_1^2 + r^2}$, this becomes
$$ Pf(x_1, t) = \omega_{n-3} \int_{\R} f(y,t) (y^2 - |x_1|^2)_+^{\frac{n-4}{2}} |y|\wrt y.$$
The estimate \eqref{l1} can thus be rewritten as
\begin{eqnarray*}
\| \int 
f(y,t) |x_1|^{\frac{n-2}{2n}} |\partial_{x_1}|^{(n-2)/2}
(y^2 - |x_1|^2)_+^{\frac{n-4}{2}} |y|\wrt y \|_{L^{2n/(n+1)}(\R \times \R)}\\
\leq C (\int |f(y,t)|^{2n/(n+1)} |y|^{n-2}\wrt y \wrt t)^{(n+1)/2n}.
\end{eqnarray*}
Freezing $t$ and setting $h(y) := f(y,t) |y|^{(n-2)(n+1)/2n}$, we thus reduce (after some algebra) to showing the one-dimensional estimate
\begin{equation}\label{k+}
\| \int 
K_+(x, y) h(y) \wrt y \|_{L^{2n/(n+1)}(\R)}
\leq C \| h \|_{L^{2n/(n+1)}(\R)},
\end{equation}
where the kernel $K_+$ is defined by
$$ K_+(x, y) := |y|^{-\frac{n-4}{2}}
(|x|/|y|)^{\frac{n-2}{2n}} |\partial_{x}|^{(n-2)/2} 
(y^2 - x^2)_+^{\frac{n-4}{2}}.$$
Similarly, to prove \eqref{l2} it will suffice to show that
\begin{equation}\label{k-}
\| \int 
K_-(x, y) h(y) \wrt y \|_{L^{2n/(n-1)}(\R)}
\leq C \| h \|_{L^{2n/(n-1)}(\R)},
\end{equation}
where the kernel $K_-$ is defined by
$$ K_-(x, y) := |y|^{-\frac{n-4}{2}}
(|x|/|y|)^{-\frac{n-2}{2n}} |\partial_{x}|^{(n-2)/2} 
(y^2 - x^2)_+^{\frac{n-4}{2}}.$$

We now estimate the kernels $K_\pm$.  First observe that $K_\pm(x,y)$ is even in both the $x$ and $y$ variables, so we may freely restrict both variables to the positive half-line $\R^+$.  Next, we observe the scaling relationship
$$ K_\pm(\lambda x, \lambda y) = \lambda^{-1} K_\pm(x, y)$$
for all $\lambda > 0$ and $x, y \in \R^+$. Thus we have
$$ K_\pm(x, y) = \frac{1}{y} K_\pm(x/y, 1).$$
It is thus of interest to estimate $K_\pm(x, 1)$.  First suppose that $x \geq 2$.  Then from the decay of the 
kernel of the pseudo-differential operator $|\partial_x|^{(n-2)/2}$ we have
\begin{equation}\label{decay}
|\partial_{x}|^{(n-2)/2} (1 - x^2)_+^{\frac{n-4}{2}} = O(|x|^{-n/2}).
\end{equation}
When $0 < x < 2$ one has to be more careful.  First observe from standard stationary phase that the Fourier transform of $|\partial_{x}|^{(n-2)/2} (1 - x^2)_+^{\frac{n-4}{2}}$ (thought of as a function on $\R$) is even, real-valued and of the form
$$ |\xi|^{(n-2)/2} \Re( C e^{2\pi i \xi}/|\xi|^{(n-2)/2} + C e^{2\pi i \xi}/|\xi|^{n/2} + O( 1/|\xi|^{(n+2)/2} )$$
for $\xi \geq 1$, where $C$ denotes various absolute constants which vary from line to line.  The error term is integrable, and so we see from inverting the Fourier transform again that
$$ |\partial_{x}|^{(n-2)/2} (1 - x^2)_+^{\frac{n-4}{2}} = 
C \delta(x^2 - 1) + C \mbox{\rm{p.v.}} \frac{1}{x^2 - 1} + C (1 - x^2)_+^0 + C + O(1),$$
where $\delta(x)$ now denotes the Dirac delta, and $\mbox{\rm{p.v.}} \frac{1}{x}$
is the distributional kernel of the Hilbert transform.  In particular, we have
\begin{equation}\label{asym}
|\partial_{x}|^{(n-2)/2} (1 - x^2)_+^{\frac{n-4}{2}} = C \delta(x - 1) + C \mbox{\rm{p.v.}} \frac{1}{x-1} + O(1),
\end{equation}
when $0 < x < 2$.

We now estimate $K_\pm(x,y)$ in the three regions $0 < x < y/2$, $y/2 < x < 2y$, and $x > 2y > 0$, in order to establish the bounds \eqref{k+}, \eqref{k-}.
First suppose we are in the region $y/2 < x < 2y$.  Then from \eqref{asym} (and a Taylor expansion) we have
$$ K_\pm(x,y) = C \delta(x-y) + C \mbox{\rm{p.v.}} \frac{1}{x-y} + O(\frac{1}{x}).$$
>From Hardy's inequality, see \cite[\S 9.8, p. 239]{HLP}\footnote{It is interesting to note that if $f(x)\ge 0$ and $f(x)=0$ for $x\le 0$ then for $x\ge 0$ $$\frac{1}{x}\int_0^xf(-t)dt \le 2\int_x^{2x}\frac{f(x-y)}{y}dy \le 2\int_{\R}\frac{f(x-y)}{y}dy.$$ Hence Hardy's inequality can be obtained directly from the $L^p$ boundedness of the Hilbert transform.},
$$ \| \frac{1}{x} \int_0^x |f(y)|\wrt y \|_p \leq C_p \| f\|_p \hbox{ for } 1 < p < \infty $$
and the $L^p$ boundedness of the Hilbert transform we thus see that the portion of $K_\pm$ in this region is acceptable.

Next, suppose we are in the region $x > 2y$.  Then from \eqref{decay} we have
$$ K_\pm(x,y) = O( \frac{1}{y} (\frac{x}{y})^{\pm (n-2)/2n - n/2} ) = O(1/x).$$
Thus this contribution is acceptable from Hardy's inequality.

Next suppose that we are in the region $x < y/2$.  Then from \eqref{asym} we have
$$ K_\pm(x,y) = O( \frac{1}{y} (\frac{x}{y})^{\pm (n-2)/2n} ) .$$
For \eqref{k+} this is again acceptable by (the adjoint of) Hardy's inequality, since the right-hand side is $O(1/y)$.  For \eqref{k-} we have to show
\begin{equation}\label{k--}
\| \int_{x < y/2} \frac{1}{y} (\frac{x}{y})^{-(n-2)/2n} |f|(y)\wrt y \|_{2n/(n-1)} \lesssim \| f \|_{2n/(n-1)}.
\end{equation}
First consider the portion of the integral where $2^{-j-1} y \leq x < 2^{-j} y$ for some integer $j > 0$.  We can estimate this portion by
$$ 2^{(n-2)j/2n} \| \int_{x < y/2^j}  \frac{1}{y}|f|(y)\wrt y \|_{2n/(n-1)},$$
which after a change of variables becomes
$$ 2^{(n-2)j/2n} 2^{-(n-1)j/2n} \|\int_{x < y}  \frac{1}{y}|f|(y)\wrt y \|_{2n/(n-1)}.$$
Summing in $j$ and using Hardy's inequality again, we obtain \eqref{k--}.

\end{proof}

\subsection*{Cylindrically symmetric functions on $\R^3 \times \R$}

We now remark that the above proof can be simplified in the special case $n=4$.  The key observation is that for functions $f(x)$ radial in $\R^3$, the three-dimensional Laplacian $\Delta_{\R^3}$ and the one-dimensional Laplacian $\Delta_\R$ are intertwined by the well-known identity
$$ \Delta_{\R^3} f(x_1 e_1) = x_1^{-1} \Delta_\R (x_1 f(x_1 e_1))$$
where $\Delta_\R = \partial_{x_1}^2$ is the Laplacian in the $x_1$ variable.
Indeed, this identity easily follows from the polar representation $\Delta_{\R^3} = \partial_r^2 + \frac{2}{r} \partial_r = r^{-1} \partial_r^2 r$ of the three-dimensional Laplacian on radial functions.

As a consequence of this identity, we see that the four-dimensional Laplacian $\Delta_{\R^3 \times \R}$ on cylindrically symmetric functions $f(x,t)$ obeys the identity
$$ \Delta_{\R^3 \times \R} f(x_1 e_1, t) = x_1^{-1} \Delta_{\R \times \R}( x_1 f(x_1 e_1, t) ).$$
In particular, the functional calculi of $\Delta_{\R^3 \times \R}$ and $\Delta_{\R \times \R}$ intertwine and we have
$$ S^\delta_1(\Delta_{\R^3 \times \R}) f(x_1 e_1, t) = x_1^{-1} S^\delta_1(\Delta_{\R \times \R})( x_1 f(x_1 e_1, t) ).$$
We need to prove \eqref{dd0}, which in polar co-ordinates (using $x_1$ as a proxy for the radial variable) becomes
$$
\int |S^\delta_1( \Delta_{\R^{n-1} \times \R} ) f(x_1 e_1, t) )|^{8/5} |x_1|^2\wrt x_1 \wrt t \leq C
\int | f(x_1 e_1, t)|^{8/5} |x_1|^2\wrt x_1 \wrt t.$$
Applying the above intertwining identity, this becomes
$$
\int |S^\delta_1( \Delta_{\R \times \R} ) g(x_1, t) )|^{8/5} |x_1|^{2/5}\wrt x_1 \wrt t \leq C
\int | g(x_1 e_1, t)|^{8/5} |x_1|^{2/5}\wrt x_1 \wrt t,$$
or in other words that $|x_1|^{1/4} S^\delta_1(\Delta_{\R \times \R}) |x_1|^{-1/4}$ is bounded on $L^{8/5}(\R \times \R)$.  But this follows from Lemma \ref{delta}.

The above argument is of course very similar to the one in the previous section, using many of the same tools.  Indeed if one ran the previous section argument for $n=4$, but distributed the derivative $|\xi_1|^2$ using $\partial_{x_1}$ instead of $|\partial_{x_1}|$, then the inequalities \eqref{l1}, \eqref{l2} would become trivial and the two arguments become essentially identical.

\end{document}